\documentclass{article}

\pdfoutput=1

\usepackage{microtype}
\usepackage{graphicx}
\usepackage{subfigure}
\usepackage{booktabs} 

\usepackage{hyperref}

\usepackage[accepted]{icml2018}

\usepackage{amsmath,amsbsy,amsgen,amscd,amssymb,amsthm,amsfonts,stmaryrd}
\usepackage{mathtools}

\usepackage{bm}

\usepackage{microtype}

\usepackage{hyperref}
\usepackage{url}

\usepackage[usenames,dvipsnames]{xcolor}
\usepackage{graphicx}

\graphicspath{{art/}}

\definecolor{dkblue}{rgb}{0,0,.5}
\definecolor{rust}{rgb}{0.5,0.1,0.1}

\hypersetup{urlcolor=rust}
\hypersetup{citecolor=dkblue}
\hypersetup{linkcolor=dkblue}

\newtheorem{theorem}{Theorem}

\theoremstyle{definition}
\newtheorem{remark}[theorem]{Remark}

\numberwithin{equation}{section}
\numberwithin{theorem}{section}

\newcommand{\R}{\mathbb{R}}

\newcommand{\trace}{\operatorname{tr}}

\newcommand{\abs}[1]{\vert #1 \vert}
\newcommand{\norm}[1]{\Vert #1 \Vert}
\newcommand{\ip}[2]{\big\langle #1, \, #2 \big\rangle}
\newcommand{\prox}[1]{\mathrm{prox}_{#1}}
\newcommand{\proj}[1]{\mathrm{proj}_{#1}}

\newcommand{\subjto}{\text{subject to}}

\DeclareFontFamily{OT1}{pzc}{}
\DeclareFontShape{OT1}{pzc}{m}{it}{<-> s * [1.200] pzcmi7t}{}
\DeclareMathAlphabet{\mathpzc}{OT1}{pzc}{m}{it}

\usepackage[nameinlink]{cleveref}
\crefname{lemma}{lemma}{lemmas}
\newcommand*{\bfrac}[2]{\genfrac{[}{]}{0pt}{}{#1}{#2}}

\icmltitlerunning{A Conditional Gradient Framework for Composite Convex Minimization with Applications to Semidefinite Programming}

\begin{document}

\twocolumn[
\icmltitle{A Conditional Gradient Framework for Composite Convex Minimization \\ with Applications to Semidefinite Programming}

\icmlsetsymbol{equal}{*}

\begin{icmlauthorlist}
\icmlauthor{Alp Yurtsever}{epfl}
\icmlauthor{Olivier Fercoq}{paris}
\icmlauthor{Francesco Locatello}{ethz,mpi}
\icmlauthor{Volkan Cevher}{epfl}
\end{icmlauthorlist}

\icmlaffiliation{epfl}{LIONS, Ecole Polytechnique F\'ed\'erale de Lausanne, Switzerland}
\icmlaffiliation{paris}{LTCI, T\'el\'ecom ParisTech, Universit\'e Paris-Saclay, France}
\icmlaffiliation{ethz}{BMI, Dept. of Computer Science, ETH Zurich, Switzerland}
\icmlaffiliation{mpi}{Empirical Inference, Max Planck Institute for Intelligent Systems, Germany}

\icmlcorrespondingauthor{Alp Yurtsever}{alp.yurtsever@epfl.ch}
\icmlkeywords{Optimization, Frank-Wolfe, Conditional gradient method, composite convex, smoothing, homotopy, quadratic assignment}

\vskip 0.3in
]

\printAffiliationsAndNotice{}  

\begin{abstract} 
We propose a conditional gradient framework for a composite convex minimization template with broad applications. 
Our approach combines smoothing and homotopy techniques under the CGM framework, and provably achieves the optimal $\mathcal{O}(1/\sqrt{k})$ convergence rate. 
We demonstrate that the same rate holds if the linear subproblems are solved approximately with additive or multiplicative error.  
In contrast with the relevant work, we are able to characterize the convergence when the non-smooth term is an indicator function.   
Specific applications of our framework include the non-smooth minimization, semidefinite programming, and minimization with linear inclusion constraints over a compact domain. 
Numerical evidence demonstrates the benefits of our framework. 
\end{abstract}

\section{Introduction}
The importance of convex optimization in machine learning has increased dramatically in the last decade due to the new theory in structured sparsity, rank minimization and statistical learning models like support vector machines. Indeed, a large  class of learning formulations can be addressed by the following composite convex minimization template: 
\begin{equation}
\label{eqn:main-template}
\min_{x \in \mathcal{X}} F(x) := f(x) + g(Ax),  
\end{equation}
where $\mathcal{X} \subset \R^n$ is compact (nonempty, bounded, closed) and its 0-dimensional faces (\textit{i.e.}, its vertices) are called \textit{atoms}. $f: \mathcal{X} \to \R \cup \{ +\infty \}$ is a smooth proper closed convex function, $A \in \R^{d \times n}$, and $g : \R^d \to \R \cup \{ +\infty \}$ is a proper closed convex function which is possibly non-smooth. 

By using the powerful proximal gradient framework, problems belonging to the template \eqref{eqn:main-template} can be solved nearly as efficiently as if they were fully smooth with fast convergence rates. By proximal ($\prox{}$) operator, we mean:
\begin{equation*}
\vspace{-0.25mm}
\prox{\phi}(v) = \arg\min_{x} ~ \phi(x) + \frac{1}{2} \norm{x-v}^2.
\vspace{-0.25mm}
\end{equation*}
These methods make use of the gradient of the smooth function $f$ along with the $\prox{}$ of the non-smooth part $g(A\,\cdot\,) + \iota_{\mathcal{X}}$, where $\iota_\mathcal{X}$ denotes the indicator function of $\mathcal{X}$, and are optimal as they match the iteration lower-bounds.  

Surprisingly, the proximal operator can impose an undesirable computational burden and even intractability on these gradient-based methods, such as the computation of a full singular value decomposition in the ambient dimension or the computation of proximal mapping for the latent group lasso \cite{Jaggi2013}. Moreover, the linear mapping $A$ often complicates the computation of the $\prox{}$ itself, and require more sophisticated splitting or primal-dual methods. 

As a result, the conditional gradient method (CGM, aka Frank-Wolfe method) has recently increased in popularity, since it requires only a linear minimization oracle ($\mathrm{lmo}$). By $\mathrm{lmo}$, we mean a resolvent of the following problem:
\begin{equation*}
\vspace{-0.25mm}
\mathrm{lmo}_{\mathcal{X}}(v) = \arg\min_{x\in\mathcal{X}} \ip{x}{v}.
\vspace{-0.25mm}
\end{equation*}
CGM features significantly reduced computational costs (\textit{e.g.}, when $\mathcal{X}$ is the spectrahedron), tractability (\textit{e.g.}, when $\mathcal{X}$ is a latent group lasso norm), and interpretability (\textit{e.g.}, they generate solutions as a combination of atoms of $\mathcal{X}$). The method is  shown in \Cref{alg:cgm} when $g(Ax)=0$: 
\vspace{-1mm}

\begin{algorithm}[h]
   \caption{CGM for smooth minimization}
   \label{alg:cgm}
\begin{algorithmic}
   \STATE {\bfseries Input:} $x_1 \in \mathcal{X}$
   \FOR{$k=1,2, \ldots, $}
   \STATE $\eta_k = 2/(k+1)$
	\STATE $s_k = \arg\min_{x \in \mathcal{X}}\ip{\nabla f(x_k)}{x}$
   	\STATE $x_{k+1} = x_k + \eta_k(s_k - x_k)$
   \ENDFOR
\end{algorithmic}
\end{algorithm}

\vspace{-1mm}
The method itself is optimal for this particular template \cite{Lan2014}. 
Unfortunately, CGM \textit{provably} cannot handle the non-smooth term $g(Ax)$ in \eqref{eqn:main-template} via its subgradients (\textit{cf.} \Cref{sec:non-smooth} for a counter example by \citet{Nesterov2017}). 

When the non-smooth part is an indicator function, one could take the intersection between $\mathcal{X}$ and the set represented by $g$. Unfortunately, even the $\mathrm{lmo}$ itself can be a difficult optimization problem depending on the structure of the domain. On many domains of interest, that we can parametrize as a composition of simple sets, linear problems are infeasible  \cite{richard2012estimation,yen2016convex}.

In this paper, we propose a CGM framework for solving the composite problem \eqref{eqn:main-template} with rigorous convergence guarantees. Our approach retains the simplicity of projection free methods, but allows to disentangle the complexity of the feasible set in order to preserve the simplicity of the $\mathrm{lmo}$. 

Our method combines the ideas of smoothing \cite{Nesterov2005} and homotopy under the CGM framework. 
\citet{Lan2014} proposes a similar approach for non-smooth problems, which is extended for the conditional gradient sliding framework in \cite{Lan2016,Lan2017}. 
Their analysis, however, is restricted by the Lipschitz continuity assumption. 
Consequently, it does not apply to the standard semidefinite programming template, or to the problems with affine inclusion constraints, limiting its applicability in machine learning (\textit{cf.} \Cref{sec:affine,sec:splitting}). 

Our work covers in particular the case where non-smooth part is the indicator function of a convex set.
Similar ideas can be found for the primal-dual subgradient method and the coordinate descent in \cite{TranDinh2017, Alacaoglu2017} via the projection onto $\mathcal{X}$. 

Our contributions can be summarized as follows: \vspace{-0.75em}
\begin{enumerate}
\itemsep-0.15em
\item[$\triangleright$] We introduce a simple, easy to implement CGM framework for solving composite problem \eqref{eqn:main-template}, and prove that it achieves the optimal $\mathcal{O}(1/\sqrt{k})$ rate. 
\item[$\triangleright$] We prove $\mathcal{O}(1/\sqrt{k})$ rate both in the objective residual and the feasibility gap, when the non-smooth term is an indicator function.
\item[$\triangleright$] We analyze the convergence of our algorithm under inexact oracles with additive and multiplicative errors. 
\item[$\triangleright$] We present key instances of our framework, including the non-smooth minimization, minimization with linear inclusion constraints, and convex splitting. 
\item[$\triangleright$] We present empirical evidence supporting our findings.
\end{enumerate}
 
\paragraph{Roadmap.}
\Cref{sec:background} recalls some basic notions and presents the preliminaries about the smoothing technique. 
In \Cref{sec:algorithm}, we present CGM for composite convex minimization along with the convergence guarantees, 
and we extend these results for inexact oracle calls in \Cref{sec:inexact}. 
We describe some important special applications of our framework in \Cref{sec:applications}. 
We provide empirical evidence supporting our theoretical findings in \Cref{sec:experiments}. 
Finally, \Cref{sec:conclusions} draws the conclusions with a discussion on the future work. 
Proofs and technical details are deferred to the appendix. 

\section{Notation \& Preliminaries}
\label{sec:background}

We use $\|\cdot\|$ to denote the Euclidean norm for vectors and the spectral norm ({a.k.a.} Schatten $\infty$-norm) for  linear mappings. 
We denote the Frobenius norm by $\|\cdot\|_F$, and the nuclear norm ({a.k.a.} Schatten 1-norm or trace norm) by $\|\cdot\|_{S_1}$. 
The notation $\ip{\cdot}{\cdot}$ refers the Euclidean or Frobenius inner product. 
The symbol $^\top$ denotes the adjoint of a linear map, and the symbol $\succcurlyeq$ denotes the semidefinite order. 
We denote the diameter of $\mathcal{X}$ by $D_{\mathcal{X}} = \max_{x_1,x_2\in\mathcal{X}}\|x_1 - x_2\|$, and the distance between a point $y$ and a set $\mathcal{K}$ by $\mathrm{dist}(y,\mathcal{K})$.

\paragraph{Lipschitz continuity.}
We say that a function $g:\R^d \to \R$ is $L$-Lipschitz continuous if it satisfies
\begin{equation*}
\abs{ g (x_1) - g(x_2)} \leq L \norm{x_1 - x_2}, \quad \forall x_1, x_2 \in \R^d. 
\end{equation*}

\paragraph{Smoothness.}
A differentiable function ${f:\mathcal{X} \to \R}$ is $L_f$-smooth if the gradient $\nabla f$ is $L_f$-Lipschitz continuous:
\begin{equation*}
\norm{\nabla f (x_1) - \nabla f(x_2)} \leq L_f \norm{x_1 - x_2}, \quad \forall x_1, x_2 \in \mathcal{X}. 
\end{equation*}

\paragraph{Fenchel conjugate \& Smoothing.}
We consider the smooth approximation of a non-smooth term $g$ obtained using the technique described by \citet{Nesterov2005} with the standard Euclidean proximity function $\tfrac{1}{2}\norm{\cdot}^2$ and a smoothness parameter $\beta > 0$
\begin{equation*}
g_\beta (z) = \max_{y \in \R^d} \ip{z}{y} - g^*(y) - \frac{\beta}{2} \norm{y}^2,
\end{equation*}
where $g^*$ denotes the Fenchel conjugate of $g$
\begin{equation*}
g^*(x) = \sup_{v} \ip{x}{v} - g(v).
\end{equation*}
Note that $g_\beta$ is convex and $\tfrac{1}{\beta}$-smooth. 

\paragraph{Solution set.}
We denote an exact solution of \eqref{eqn:main-template} by $x^\star$, and the set of all solutions by $\mathcal{X}^\star$.
Throughout the paper, we assume that the solution set $\mathcal{X}^\star$ is nonempty. 

Given an accuracy level $\epsilon > 0$, we call a point $x \in \mathcal{X}$ as an $\epsilon$-solution of \eqref{eqn:main-template} if
\begin{equation}\label{eqn:eps-solution-1}
f(x) + g(Ax) - f^\star - g^\star \leq \epsilon,
\end{equation}
where we use the notation $f^\star = f(x^\star)$ and $g^\star = g(Ax^\star)$.

When $g$ is the indicator function of a set $\mathcal{K}$, condition \eqref{eqn:eps-solution-1} is not well-defined for infeasible points. 
Hence, we refine our definition, and call a point $x \in \mathcal{X}$ as an $\epsilon$-solution if
\begin{equation*}
f(x) - f^\star \leq \epsilon, \quad \text{and} \quad \mathrm{dist}(Ax, \mathcal{K}) \leq \epsilon.
\end{equation*}
Here, we call $f(x) - f^\star$ as the objective residual and $\mathrm{dist}(Ax, \mathcal{K})$ as the feasibility gap. 
We use the same $\epsilon$ for the objective residual and the feasibility gap, since the distinct choices can be handled by scaling $f$. 

\paragraph{Lagrange saddle point.} Suppose that $g$ is the indicator function of a convex set $\mathcal{K}$, and 
denote the Lagrangian of problem \eqref{eqn:main-template} by $\mathcal{L}(x,y)$:
$$\mathcal{L}(x,y) := f(x) + \ip{y}{Ax} - g^\ast(y).$$
We can formulate primal and dual problems as follows:
$$\underbrace{\max_{y \in \R^d}  \min_{ x \in \mathcal{X}} \mathcal{L}(x,y)}_{\text{dual}} \leq \underbrace{\min_{ x \in \mathcal{X}}\max_{y \in \R^d} \mathcal{L}(x,y)}_{\text{primal}}.$$
We assume that the strong duality holds, \textit{i.e.}, the relation above holds with equality. 
Slater's condition is a sufficient condition for strong duality. 
By Slater's condition, we mean 
$$\mathrm{relint}(\mathcal{X} \times \mathcal{K}) ~\cap~ \{(x,r) \in \R^n \times \R^d : Ax = r \} \neq \emptyset,$$
where $\mathrm{relint}$ stands for the relative interior. 

Throughout, $y^\star$ denotes a solution of the dual problem. 

\section{Algorithm \& Convergence}
\label{sec:algorithm}

Our method is based on the simple idea of combining smoothing and homotopy. 
In our problem template, the objective function $F$ is non-smooth. 
We define the smooth approximation of $F$ with smoothness parameter $\beta > 0$ as 
\begin{equation*}
F_\beta(x) = f(x) + g_{\beta}(Ax).
\end{equation*}
Note that $F_\beta$ is $(L_f + \norm{A}^2/\beta)$-smooth. 

The algorithm takes a conditional gradient step with respect to the smooth approximation $F_{\beta_k}$ at iteration $k$, where $\beta_k$ is gradually decreased towards $0$. 
Let us denote by $y^\ast_{\beta_k}$
\begin{align*}
y^\ast_{\beta_k} (Ax)
& = 
\arg\max_{y \in \R^d} \ip{Ax}{y} - g^*(y) - \frac{\beta_k}{2} \norm{y}^2 \\
& = \prox{\beta_k^{-1} g^\ast} (\beta_k^{-1} Ax) \\
& = \frac{1}{\beta_k} \big( Ax - \prox{\beta_k g} (Ax) \big),
\end{align*}
where the last equality is due to the Moreau decomposition. 
Then, we can compute the gradient of $F_{\beta_k}$ as
\begin{align*}
\nabla F_{\beta_k}(x) 
& = 
\nabla f(x) + A^\top y^\ast_{\beta_k}(Ax) \\
& = 
\nabla f(x) + \frac{1}{\beta_k} A^\top \big( Ax - \prox{\beta_k g} (Ax) \big). 
\end{align*}

Based on this formulation, we present our CGM framework for composite convex minimization template \eqref{eqn:main-template} in \Cref{alg:composite}. 
The choice of $\beta_k$ comes from the convergence analysis, which can be found in the supplements. 

\begin{algorithm}[h]
   \caption{CGM for composite problems}
   \label{alg:composite}
\begin{algorithmic}
   \STATE {\bfseries Input:} $x_1 \in \mathcal{X}, ~ \beta_0 > 0$
   \FOR{$k=1,2, \ldots, $}
   \STATE $\eta_k = 2/(k+1)$,~~~and~~~$\beta_k =\beta_0 / \sqrt{k+1}$
   \STATE $v_k = \beta_k \nabla f(x_k) + A^\top \big(A x_k -  \prox{\beta_k g}(Ax_k)\big)$
	\STATE $s_k = \arg\min_{x \in \mathcal{X}}\ip{v_k}{x}$
   	\STATE $x_{k+1} = x_k + \eta_k(s_k - x_k)$
   \ENDFOR
\end{algorithmic}
\end{algorithm}

\begin{theorem}\label{thm:composite}
The sequence $x_k$ generated by \Cref{alg:composite} satisfies the following bound:
\begin{align*}
F_{\beta_k}(x_{k+1}) - F^\star 
\leq 2 D_{\mathcal{X}}^2 \left( \frac{ L_f}{k+1} + \frac{\norm{A}^2}{\beta_0 \sqrt{k+1}} \right).
\end{align*}
\end{theorem}

\Cref{thm:composite} does not directly certify the convergence of $x_k$ to a solution, since the bound is on the smoothed gap $F_{\beta_{k}}(x_k) - F^\star$. 
To relate $F_{\beta_{k}}$ back to $F$, one usually assumes Lipschitz continuity. 
This well known perspective leads us to \Cref{thm:lipschitz}, which is a direct extension of Theorem 4 of \cite{Lan2014} for composite functions. 

\begin{theorem}\label{thm:lipschitz}
Assume that $g: \R^d \to \R$ is $L_g$-Lipschitz continuous. 
Then, the sequence $x_k$ generated by \Cref{alg:composite} satisfies the following convergence bound:
\begin{align*}
F(x_k) - F^\star 
\leq 2 D_{\mathcal{X}}^2 \left( \frac{ L_f}{k} + \frac{\norm{A}^2}{\beta_0 \sqrt{k}} \right) + \frac{\beta_0 L_g^2}{2\sqrt{k}}.
\end{align*}
Furthermore, if the constants $D_{\mathcal{X}}$, $\norm{A}$ and $L_g$ are known or easy to approximate, we can choose $\beta_0 = 2D_{\mathcal{X}}\norm{A}/L_g$ to get the following convergence rate:
\begin{align*}
F(x_k) - F^\star 
\leq   \frac{ 2 D_{\mathcal{X}}^2 L_f}{k} + \frac{2D_{\mathcal{X}}\norm{A}L_g}{\sqrt{k}}.
\end{align*}
\end{theorem}

Lipschitz continuity assumption in \Cref{thm:lipschitz} leaves many important applications out ({\em cf.} \Cref{sec:affine,sec:splitting}). 
In \Cref{thm:indicator}, we take a step further and characterize the convergence when the non-smooth part is an indicator function. 
Note that $x_k$ is not guaranteed to be a feasible point, since the condition $Ax_k \in \mathcal{K}$ is not guaranteed, 
but it converges to the feasible set, and it becomes feasible asymptotically. 

\begin{theorem}\label{thm:indicator}
Assume that $g: \R^d \to \R$ is the indicator function of a simple convex set $\mathcal{K}$. 
Then, the sequence $x_k$ generated by \Cref{alg:composite} satisfies:
\begin{align*}
f(x_{k}) - f^\star & \geq -\norm{y^\star} ~ \mathrm{dist}(Ax_k,\mathcal{K}) \\
f(x_{k}) - f^\star & \leq 2 D_{\mathcal{X}}^2 \left( \frac{ L_f}{k} + \frac{\norm{A}^2}{\beta_0 \sqrt{k}} \right) \\
\mathrm{dist}(Ax_k,\mathcal{K}) & \leq  \frac{2 \beta_0}{\sqrt{k}} \left( \norm{y^\star} + D_{\mathcal{X}}\sqrt{\frac{C_0}{\beta_0}} \right)
\end{align*}
where $C_0 = L_f + \|A\|^2 / \beta_{0}$.
\end{theorem}

\begin{remark}
Similar to the classical CGM, we can consider variants of \Cref{alg:composite} with line-search (which replaces the step size by $\eta_k = \min_{\eta \in [0,1]} F_{\beta_k} (x_k + \eta(s_k - x_k))$), and fully corrective updates (which replaces the last step by $x_{k+1} = \arg\min_{x \in \textrm{conv}(s_1,\ldots,s_k)} F_{\beta_k}(x)$). 
All results presented in this paper remain valid for these variants.  
\end{remark}

\section{Convergence with Inexact Oracles} 
\label{sec:inexact}

Finding an exact solution of the $\mathrm{lmo}$ can be expensive in practice, especially when it involves a matrix factorization as in the SDP examples.
On the other hand, approximate solutions can be much more efficient. 

Different notions of inexact $\mathrm{lmo}$ are already explored in CGM and greedy optimization frameworks, \textit{cf.}\cite{lacoste2012block,locatello2017unified,locatello2017greedy}. We revisit the notions of additive and multiplicative errors which we adapt here for our setting.

\subsection{Inexact Oracle with Additive Error}
At iteration $k$, for the given direction $v_k$, we assume that the approximate $\mathrm{lmo}$ returns an element $\tilde{s}_k\in\mathcal{X}$ such that:
\begin{equation}\label{eqn:inexact_LMO}
\ip{v_k}{\tilde{s}_k} \leq \ip{v_k}{s_k} + \delta \frac{\eta_k}{2}D_{\mathcal{X}}^2\left(L_f + \frac{\|A\|^2}{\beta_k}\right)
\end{equation}
for some $\delta > 0$, where $s_k$ denotes the exact solution of the $\mathrm{lmo}$. 
Note that as in ~\cite{Jaggi2013}, we require the accuracy of $\mathrm{lmo}$ to increase as the algorithm progresses. 

Replacing the exact $\mathrm{lmo}$ with the approximate oracles of the form \eqref{eqn:inexact_LMO} in \Cref{alg:composite}, we get the convergence guarantees in \Cref{thm:composite-inexact,thm:lipschitz-inexact,thm:indicator-inexact}. 

\begin{theorem}\label{thm:composite-inexact}
The sequence $x_k$ generated by \Cref{alg:composite} with approximate $\mathrm{lmo}$ \eqref{eqn:inexact_LMO} satisfies:
\begin{align*}
F_{\beta_k}(x_{k+1}) - F^\star 
\leq 2 D_{\mathcal{X}}^2 \left( \frac{ L_f}{k+1} + \frac{\norm{A}^2}{\beta_0 \sqrt{k+1}} \right)(1+\delta).
\end{align*}
\end{theorem}

\begin{theorem}\label{thm:lipschitz-inexact}
Assume that $g$ is $L_g$-Lipschitz continuous. 
Then, the sequence $x_k$ generated by \Cref{alg:composite} with approximate $\mathrm{lmo}$ \eqref{eqn:inexact_LMO} satisfies:
\begin{align*}
F(x_k) - F^\star 
\leq 2 D_{\mathcal{X}}^2 \left( \frac{ L_f}{k} + \frac{\norm{A}^2}{\beta_0 \sqrt{k}} \right)(1+\delta) + \frac{\beta_0 L_g^2}{2 \sqrt{k}}.
\end{align*}
We can optimize $\beta_0$ from this bound if $\delta$ is known. 
\end{theorem}

\begin{theorem}\label{thm:indicator-inexact}
Assume that $g$ is the indicator function of a simple convex set $\mathcal{K}$. 
Then, the sequence $x_k$ generated by \Cref{alg:composite} with approximate $\mathrm{lmo}$ \eqref{eqn:inexact_LMO} satisfies:
\begin{align*}
f(x_{k}) - f^\star & \geq -\norm{y^\star} \mathrm{dist}(Ax_{k} , \mathcal{K}) \\
f(x_{k}) - f^\star & \leq 2 D_{\mathcal{X}}^2 \left( \frac{ L_f}{k} + \frac{\norm{A}^2}{\beta_0 \sqrt{k}} \right)(1+\delta) \\
\mathrm{dist}(Ax_{k} , \mathcal{K}) & \leq  \frac{2 \beta_0}{\sqrt{k}} \left( \norm{y^\star} + D_{\mathcal{X}}\sqrt{\frac{C_0}{\beta_0} (1+\delta)} \right)
\end{align*}
where $C_0 = L_f + \|A\|^2 / \beta_{0}$.
\end{theorem}

\subsection{Inexact Oracle with Multiplicative Error}

We consider the multiplicative inexact oracle:
\begin{align}\label{eqn:inexact_LMO_mult}
\ip{v_k}{ \tilde{s}_k - x_k} \leq \delta\ip{v_k}{s_k-x_k}
\end{align}
where $\delta\in(0,1]$. 
Replacing the exact $\mathrm{lmo}$ with the approximate oracles of the form \eqref{eqn:inexact_LMO_mult} in \Cref{alg:composite}, we get the convergence guarantees in \Cref{thm:composite-inexact-mult,thm:lipschitz-inexact-mult,thm:indicator-inexact-mult}. 

\begin{theorem}\label{thm:composite-inexact-mult}
The sequence $x_k$ generated by \Cref{alg:composite} with approximate $\mathrm{lmo}$ of the form \eqref{eqn:inexact_LMO_mult}, and modifying $\eta_k = \frac{2}{\delta (k-1)+2}$ and $\beta_k = \frac{\beta_0}{\sqrt{\delta k+1}}$ satisfies:
\begin{align*}
F_{\beta_k}(x_{k+1}) - F^\star 
\leq 
\frac{2 }{\delta} \left( 
\frac{D_{\mathcal{X}}^2 L_f + \delta\mathcal{E}}{\delta k+2} + \frac{D_{\mathcal{X}}^2\norm{A}^2}{\beta_0\sqrt{\delta k+2}} \right)
\end{align*}
\end{theorem}
where $\mathcal{E} = F(x_1) - F^\star$.

\begin{theorem}\label{thm:lipschitz-inexact-mult}
Assume that $g$ is $L_g$-Lipschitz continuous. 
Then, the sequence $x_k$ generated by \Cref{alg:composite} with approximate $\mathrm{lmo}$ \eqref{eqn:inexact_LMO_mult}, and modifying $\eta_k = \frac{2}{\delta (k-1)+2}$ and $\beta_k = \frac{\beta_0}{\sqrt{\delta k+1}}$ satisfies:\\
\begin{align*}
F(x_k) \!-\! F^\star 
\leq \frac{2 }{\delta} \! \left( 
\frac{D_{\mathcal{X}}^2 L_f \!+\! \delta\mathcal{E}}{\delta k\!+\! 1} \!+\! \frac{D_{\mathcal{X}}^2\norm{A}^2}{\beta_0\sqrt{\delta k \!+\! 1}} \right) \!+\! \frac{\beta_0 L_g^2}{2 \sqrt{\delta k \!+\! 1}},
\end{align*}
where $\mathcal{E} = F(x_1) - F^\star$. We can optimize $\beta_0$ from this bound if $\delta$ is known.
\end{theorem}

\begin{theorem}\label{thm:indicator-inexact-mult}
Assume that $g$ is the indicator function of a simple convex set $\mathcal{K}$. 
Then, the sequence $x_k$ generated by \Cref{alg:composite} with approximate $\mathrm{lmo}$ \eqref{eqn:inexact_LMO_mult}, and modifying $\eta_k = \frac{2}{\delta (k-1)+2}$ and $\beta_k = \frac{\beta_0}{\sqrt{\delta k+1}}$ satisfies:
\begin{align*}
f(x_{k}) - f^\star & \geq -\norm{y^\star} \mathrm{dist}(Ax_{k} , \mathcal{K}) \\
f(x_{k}) - f^\star & \leq \frac{2 }{\delta} \left(\frac{D_{\mathcal{X}}^2 L_f + \delta\mathcal{E}}{\delta k+1} + \frac{D_{\mathcal{X}}^2\norm{A}^2}{\beta_0\sqrt{\delta k+1}} \right) \\
\mathrm{dist}(Ax_{k} , \mathcal{K}) & \leq  \frac{2 \beta_0}{\sqrt{\delta k + 1}} \left( \norm{y^\star} + \sqrt{\frac{D_{\mathcal{X}}^2 C_0 + \delta\mathcal{E}}{\beta_0 \delta} } \right)
\end{align*}
where $\mathcal{E} = F(x_1) - F^\star$ and $C_0 = L_f + \|A\|^2 / \beta_{0}$.
\end{theorem}

\section{Applications \& Related Work}
\label{sec:applications}

CGM is proposed for the first time in the seminal work of \citet{FrankWolfe1956} for solving smooth convex optimization on a polytope. 
It is then progressively generalized for more general settings in \cite{Levitin1966, Dunn1978, Dunn1979, Dunn1980}. 
Nevertheless, with the introduction of the fast gradient methods with $\mathcal{O}(1/k^2)$ rate by \citet{Nesterov1983},  the development of CGM-type methods entered into a stagnation period. 

The recent developments in machine learning applications with vast data brought the scalability of the first order methods under scrutiny. 
As a result, there has been a renewed interest in CGM in the last decade. 
We compare our framework with the recent developments of CGM literature in different camps of problem templates below. 

\subsection{Smooth Problems}
CGM is extended for the smooth convex minimization over the simplex by \citet{Clarkson2010}, for the spactrahedron by \citet{Hazan2008}, and for an arbitrary compact convex set by \citet{Jaggi2013}. 
Online, stochastic and block coordinate variants of CGM are introduced by \citet{Hazan2012}, \citet{Hazan2016} and \citet{lacoste2012block} respectively. 

When applied to smooth problems, \Cref{alg:composite} is equivalent to the classical CGM, and \Cref{thm:lipschitz} recovers the known optimal $\mathcal{O}(1/k)$ convergence rate. 
We refer to \cite{Jaggi2013} for a review of applications of this template. 

It needs to be mentioned that \citet{Nesterov2017} relaxes the smoothness assumption showing that CGM converges for weakly-smooth objectives (\textit{i.e.}, with H\"older continuous gradients of order $\nu \in (0,1]$). 

\subsection{Regularized Problems}
CGM for composite problems is considered recently by \citet{Nesterov2017} and \citet{Xu2017}. 
A similar but slightly different template, where $\mathcal{X}$ and $g$ are assumed to be a closed convex cone and a norm respectively, is also studied by \citet{Harchaoui2015}. 
However, these works are based on the resolvents of a modified oracle,
\begin{equation*}
\arg\min_{x\in\mathcal{X}} \ip{x}{v} + g(Ax),
\end{equation*}
which can be expensive, unless $\mathcal{X} \equiv \mathbb{R}^n$, or $g = 0$.

\Cref{alg:composite} applies to the problem template \eqref{eqn:main-template} by leveraging $\mathrm{prox}$ of the regularizer and $\mathrm{lmo}$ of the domain independently. 
This allows us to consider additional sparsity, group sparsity and structured sparsity promoting regularizations, elastic-net regularization, total variation regularization and many others under the CGM framework.

Semi-proximal mirror-prox proposed by \citet{He2015} is also based on the smoothing technique, yet the motivation is fundamentally different. 
This method considers the regularizers for which the $\prox{}$ is difficult to compute, but can be approximated via CGM. 

\subsection{Non-Smooth Problems} \label{sec:non-smooth}
Template \eqref{eqn:main-template} covers the non-smooth convex minimization template as a special case:
\begin{equation}
\label{eqn:template-nonsmooth}
\min_{x \in \mathcal{X}} g(Ax). 
\end{equation}
Unfortunately, the classical CGM (\Cref{alg:cgm}) cannot handle the non-smooth minimization, as shown by \citet{Nesterov2017} with the following counter-example. 

\textit{Example.} Let $\mathcal{X}$ be the unit Euclidean norm ball in $\R^2$, and $g(x) = \max\{x_{(1)},x_{(2)}\}$. 
Clearly, $x^\star = [-\frac{1}{\sqrt{2}}, -\frac{1}{\sqrt{2}}]^\top$.
Choose an initial point $x_0 \neq x^\star$. 
We can use an oracle that returns a subgradient $\nabla f(x) \in \bfrac{1}{0} , \bfrac{0}{1} \}$ at any point $x \in \mathcal{X}$. 
Therefore, $\mathrm{lmo}$ returns $\bfrac{\text{-}1}{0} $ or $\bfrac{0}{\text{-}1} $ at each iteration, and $x_k$ belongs to the convex hull of $\{x_0, \bfrac{\text{-}1}{0}, \bfrac{0}{\text{-}1} \}$ which does not contain the solution.

Our framework escapes such issues by leveraging $\mathrm{prox}$ of the objective function $g$ (see appendix for numerical illustration). 
In this pathological example, $\mathrm{prox}_g$ corresponds to the projection onto the simplex. 
Often times the cost of $\mathrm{prox}_g$ is negligible in comparison to the cost of $\mathrm{lmo}_\mathcal{X}$ (\textit{cf.} \Cref{sec:RobustPCA} for a robust PCA example).

Assume that $g:\R^d \to \R$ is $L_g$-Lipschitz continuous. 
As a consequence of \Cref{thm:lipschitz}, \Cref{alg:composite} for solving \eqref{eqn:template-nonsmooth} by choosing $\beta_0 = 2D_{\mathcal{X}}\norm{A}/L_g$ satisfies
\begin{align*}
g(A x_k) - g^\star 
\leq  \frac{2D_{\mathcal{X}}\norm{A}L_g}{\sqrt{k}}.
\end{align*}
We recover the method proposed by \citet{Lan2014} in this specific setting. 
\citet{Lan2014} shows that  this rate is optimal for algorithms approximating the solution of \eqref{eqn:template-nonsmooth} as a convex combination of $\mathrm{lmo}$ outputs. 

We extend the analysis in this setting for inexact oracles. 
In contrast to the smooth case, where the additive error should decrease by $\mathcal{O}({1}/{k})$ rate, 
definition \eqref{eqn:inexact_LMO} implies that we can preserve the convergence rate in the non-smooth case if the additive error is $\mathcal{O}({1}/{\sqrt{k}})$.

\subsection{Minimax Problems}
We consider the minimax problems of the following form:
\begin{equation*}
\label{eqn:saddle-point}
\min_{x \in \mathcal{X}} \max_{y \in \mathcal{Y}} \mathcal{L}(Ax,y)  
\end{equation*}
where $\mathcal{L}$ is a smooth convex-concave function, \textit{i.e.}, $\mathcal{L}(\cdot, y)$ is convex $\forall y \in \mathcal{Y}$ and $\mathcal{L}(Ax, \cdot)$ is concave $\forall x \in \mathcal{X}$. 
Note that this formulation is a special instance of \eqref{eqn:template-nonsmooth} with $g(Ax) = \max_{y \in \mathcal{Y}} \mathcal{L}(Ax,y)$.
Consequently, we can apply \Cref{alg:composite} if $\prox{g}$ is tractable. 

When $\mathcal{Y}$ admits an efficient projection oracle, $\prox{g}$ is also efficient for bilinear saddle point problems $\mathcal{L}(Ax,y) = \ip{Ax}{y}$. 
By Moreau decomposition, we have
\begin{equation*}
\prox{g}(Ax_k) = Ax_k - \proj{\mathcal{Y}}(Ax_k),
\end{equation*}
hence $v_k$ takes the form 
\begin{equation*}
v_k = \beta_k \nabla f(x_k) + A^\top \proj{\mathcal{Y}}(Ax_k).
\end{equation*}

\citet{Gidel2017} proposes a CGM variant for the smooth convex-concave saddle point problems. 
This method processes both $x$ and $y$ via the $\mathrm{lmo}$, and hence it also requires $\mathcal{Y}$ to be bounded. 
Our method, on the other hand, is more suitable when $\proj{\mathcal{Y}}$ is easy. 

Bilinear saddle point problem covers 
the maximum margin estimation of structured output models \cite{Taskar2006} and  
minimax games \cite{VonNeumann1944}. 
In particular, it also covers an important semidefinite programming formulation \cite{Garber2016}, where $\mathcal{X}$ is a spactrahedron and $\mathcal{Y}$ is the simplex. 
Our framework fits perfectly here since the projection onto the simplex can be computed efficiently. 
We defer the extension of our framework with the entropy Bregman smoothing for future. 

Note that CGM is  applicable also for the variational inequality problems beyond \eqref{eqn:main-template}, see \cite{Hammond1984}, \cite{Juditsky2016} and \cite{Cox2017}. 

\subsection{Problems with Affine Constraints}
\label{sec:affine}
\Cref{alg:composite} also applies to smooth convex minimization problems with affine constraints over a convex compact set:
\begin{equation}
\label{eqn:template-affine-constraint}
\min_{x \in \mathcal{X}} f(x) \quad \subjto \quad Ax = b, 
\end{equation}
by setting $g(Ax)$ in \eqref{eqn:main-template} as indicator function of set $\{b\}$, where $b \in \R^d$ is a known vector. 

Since the $\prox{}$ operator of the indicator function of a convex set is the projection, $v_k$ in \Cref{alg:composite} becomes
\begin{equation*}
v_k = \beta_k \nabla f(x_k) + A^\top(Ax_k-b).
\end{equation*}

\citet{Gidel2018} recently proposed a CGM framework via augmented Lagrangian for $Ax=0$ constraint. 
Their method achieves $\mathcal{O}(1/\sqrt{k})$ rate in feasibility gap, and $\mathcal{O}(1/k)$ rate in augmented Lagrangian residual. 
Note however, these alone do not imply convergence in the objective residual. 
Moreover, the step-size of the proposed method depends on the unknown error bound constant ({\em cf.} \cite{Bolte2017} for more details on error bounds).

\citet{Liu2018} introduces an inexact augmented Lagrangian method, where the subproblems are approximately solved by CGM up to an accuracy. 
This method produces an $\epsilon$-solution in $\mathcal{O}(1/\epsilon^2)$ iterations (see Corollary~4.4 in \cite{Liu2018}), 
but it calls $\mathrm{lmo}$ an unknown number of times bounded by $\mathcal{O}(k^2)$ at $k^{\text{th}}$ iteration.

Another relevant approach here is the universal primal-dual gradient method (UPD) proposed by \citet{Yurtsever2015}. 
UPD takes advantage of Fenchel-type oracles, which can be thought as a generalization of $\mathrm{lmo}$. 

Unfortunately, UPD iterations explicitly depend on the target accuracy level $\epsilon$, which is difficult to tune without a rough knowledge of the optimal value. 
Moreover, the method converges only up to $\epsilon$-suboptimality. 
There is no known analysis with inexact oracle calls for UPD, and errors in function evaluation can cause the algorithm to get stuck in the line-search procedure. 

We can generalize \eqref{eqn:template-affine-constraint} for the problems with affine inclusion constraints:
\begin{equation}
\label{eqn:template-inclusion-constraint}
\min_{x \in \mathcal{X}} f(x) \quad \subjto \quad Ax - b \in \mathcal{K}, 
\end{equation}
where $\mathcal{K}$ is a simple closed convex set. 
In this case, $v_k$ takes the following form:
\begin{equation*}
v_k = \beta_k \nabla f(x_k) + A^\top \big(Ax_k-b-\mathrm{proj}_{\mathcal{K}}(Ax_k - b)\big).
\end{equation*}

We implicitly assume that $\proj{\mathcal{K}}$ is tractable. 
We can use a splitting framework when it is computationally more advantageous to use $\mathrm{lmo}_{\mathcal{K}}$ instead (\textit{cf.} \Cref{sec:splitting}). 

This template covers the standard semidefinite programming in particular. 
Applications include 
clustering \cite{Peng2007}, 
optimal power-flow \cite{Lavaei2012}, 
sparse PCA \cite{dAspremont2004}, 
kernel learning \cite{Lanckriet2004}, 
blind deconvolution \cite{Ahmed2014}, 
community detection \cite{Bandeira2016}, 
\textit{etc.} 
Besides machine learning applications, this formulation has a crucial role in the convex relaxation of combinatorial problems. 

A significant example is the problems over the doubly nonnegative cone (\textit{i.e.}, the intersection of the positive semidefinite cone and the positive orthant) with a bounded trace norm \cite{Yoshise2010}. 
Note that the $\mathrm{lmo}$ over this domain can be costly since the $\mathrm{lmo}$ can require full dimensional updates \cite{Jester1996, locatello2017greedy}.

Our framework can handle these problems ensuring the positive semidefiniteness by $\mathrm{lmo}_\mathcal{X}$, and can still ensure the convergence to the first orthant via $\proj{\mathcal{K}}$. 
To the best of our knowledge, our framework is the first CGM extension that can handle affine constraints.  

\subsection{Minimization via Splitting} \label{sec:splitting}
We can take advantage of splitting since we can handle affine constraints. 
This lets us to disentangle the complexity of the constraints. 
Consider the following optimization template:
\begin{equation*}
\begin{aligned}
& \min_{\substack{x \in \mathcal{X}_1 \cap \mathcal{X}_2} }& & f(x) + g (A x) \\
& \subjto & & B x - b \in \mathcal{K} ,  \quad C x - c \in \mathcal{S}  
\end{aligned}
\end{equation*}
where $\mathcal{X}_1$ and $ \mathcal{X}_2 \subset \R^n$ are two convex compact sets, $A,B,C$ are known matrices and $b,c$ are given vectors. 

\begin{figure*}[ht]
\centering
\includegraphics[width=0.95\textwidth]{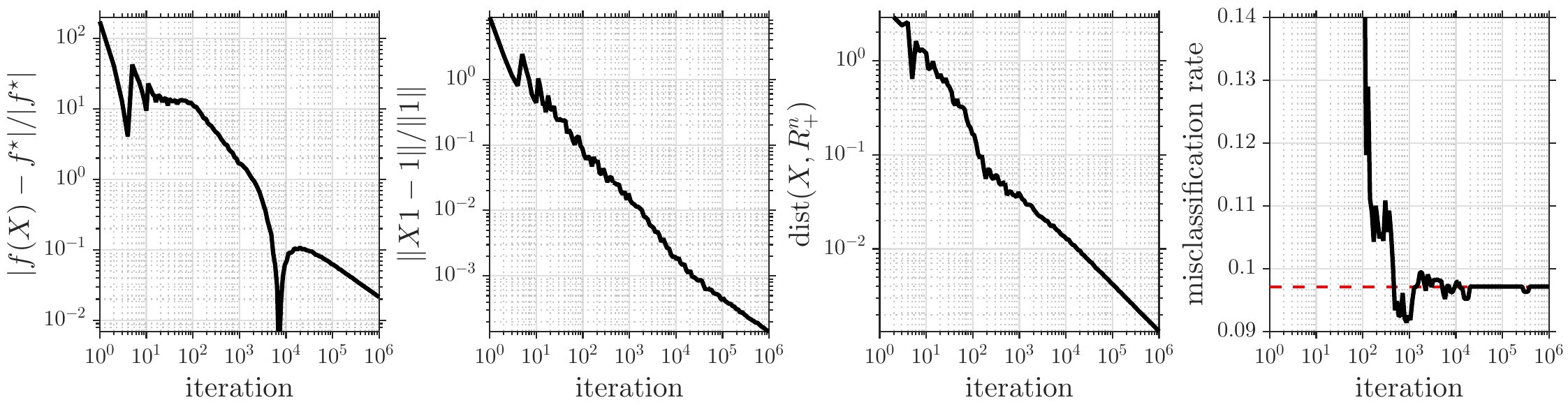}
        \vspace{-.5em}
        \caption{Clustering MNIST dataset: Convergence of our framework in function value and the feasibility gap. Red dashed line on the misclassification plot represents the value reported by \citet{Mixon2017}.}
        \label{fig:clustering}
\end{figure*}

Suppose that \\
$\cdot$~~$\mathrm{lmo}_{\mathcal{X}_1}$ and $\mathrm{lmo}_{\mathcal{X}_2}$ are easy to compute, but not $\mathrm{lmo}_{\mathcal{X}_1 \cap \mathcal{X}_2}$ \\
$\cdot$~~$\prox{g}$ is easy to compute \\
$\cdot$~~$\mathcal{K}$ is a simple convex set and $\proj{\mathcal{K}}$ is efficient \\
$\cdot$~~$\mathcal{S}$ is a convex compact set with an efficient $\mathrm{lmo}$.\\

We can reformulate this problem introducing slack variables $\xi \in \mathcal{X}_2$ and $\psi \in \mathcal{S}$ as follows:
\begin{equation*}
\begin{aligned}
& \min_{\substack{x \in \mathcal{X}_1 ~ \xi \in \mathcal{X}_2 \\ \psi \in \mathcal{S}}  }& & f(x) + g (A x) \\
& \subjto & & B x - b \in \mathcal{K}, \quad C x - c = \psi  , \quad x = \xi.
\end{aligned}
\end{equation*}
This formulation is in the form of \eqref{eqn:template-inclusion-constraint} with respect to the variable $(x,\xi,\psi) \in \mathcal{X}_1 \times \mathcal{X}_2 \times \mathcal{S}$. 
It is easy to verify that \Cref{alg:composite} leverages $\mathrm{lmo}_{\mathcal{X}_1}$, $\mathrm{lmo}_{\mathcal{X}_2}$, $\mathrm{lmo}_{\mathcal{S}}$, $\prox{g}$ and $ \proj{\mathcal{K}}$ separately. 
This approach can be generalized for an arbitrary finite number of non-smooth terms in a straightforward way.

\section{Numerical Experiments}
\label{sec:experiments}

This section presents numerical experiments supporting our theoretical findings in clustering and robust PCA examples. 
The non-smooth parts in the chosen examples consist of indicator functions, for which the dual domain is unbounded. 
Hence, to the best of our knowledge, other CGM variants in the literature are not applicable. 

\subsection{Clustering the MNIST dataset}

We consider the model-free $k$-means clustering based on the semidefinite relaxation of \citet{Peng2007}:
\begin{equation}
\label{eqn:example-clustering}
\min_{X \in \mathcal{X}} \ip{D}{X} \quad \subjto \quad \underbrace{X1 = 1, \quad X \geq 0}_{g},  
\vspace{-0.5em}
\end{equation}
where $\mathcal{X} = \{ X \in \R^{n \times n} : X \succcurlyeq 0, ~\trace(X) \leq \rho \}$ is the set of positive semidefinite matrices with a bounded trace norm, and $D \in \R^{n\times n}$ is the Euclidean distance matrix.

We use the setup described and published online by \citet{Mixon2017}, which can be briefly described as follows:
First, meaningful features from MNIST dataset \cite{MNIST}, which consists of $28 \times 28$ grayscale images that can be stacked as $784 \times 1$ vectors, are extracted using a one-layer neural network.  
This gives us a weight matrix $W \in \R^{784 \times 10}$ and a bias vector $b \in \R^{10}$. 
Then, the trained neural network is applied to the first $1000$ elements of the test set, which gives the probability vectors for these $1000$ test points, where each entry represents the probability of being each digit.

\citet{Mixon2017} runs a relax-and-round algorithm which solves \eqref{eqn:example-clustering} by SDPNAL+ \cite{Yang2015} followed by a rounding scheme (see Section~5 of \cite{Mixon2017} for details), and compares the results against MATLAB's built-in $k$-means++ implementation. 
Relax-and-round method is reported to achieve a misclassification rate of $0.0971$. 
This rate matches with the all-time best rate for $k$-means++ after $100$ different runs with random initializations. 

\begin{figure*}[ht]
\centering
\includegraphics[width=0.95\textwidth]{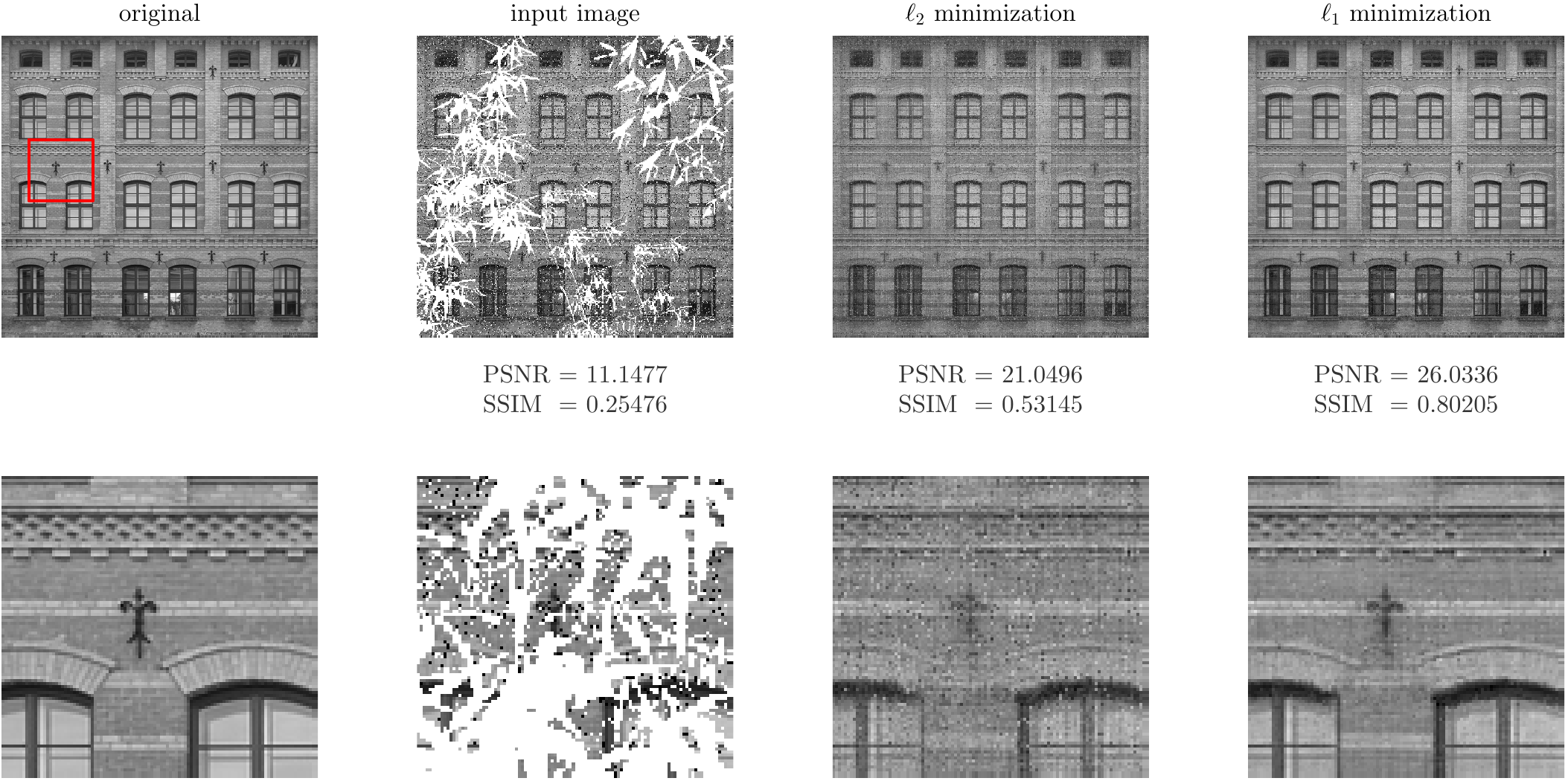}
        \vspace{-0.25em}
        \caption{Image inpainting from noisy test image $(493 \times 517)$: Robust PCA recovers a better approximation with $5$dB higher PSNR. }
        \label{fig:inpainting}
\end{figure*}

For this experiment, we solve \eqref{eqn:example-clustering} by using \Cref{alg:composite}. 
Then, we cluster data using the same rounding scheme as \cite{Mixon2017}. 
We initialize our method from zeros, and we choose $\beta_0 = 1$. 
We implement $\mathrm{lmo}$ using the built-in MATLAB function \texttt{eigs} with tolerance parameter $10^{-9}$. 

We present the results of this experiment in \Cref{fig:clustering}. 
We observe empirical $\mathcal{O}(1/\sqrt{k})$ rate both in the objective residual and the feasibility gap. 
Surprisingly, the method achieves the best test error around $1000$ iterations achieveing the misclassification rate of $0.0914$. 
This improves the value reported by \citet{Mixon2017} by $5.8 \%$. 

This example suggests that the slow convergence rate is not a major problem in many machine learning problems, since a low accuracy solution can generalize as well as the optimal point in terms of the test error, if not better.

\subsection{Robust PCA}\label{sec:RobustPCA}

Suppose that we are given a large matrix that can be decomposed as the summation of a low-rank and a sparse (in some representation) matrix. 
Robust PCA aims to recover these components accurately. 
Robust PCA has many applications in machine learning and data science, such as collaborative filtering, system identification, genotype imputation, etc. 
Here, we focus on an image decomposition problem so that we can visualize the decomposition error results.

Our setting is similar to the setup described in \cite{Zeng2018}. 
We consider a scaled grayscale photograph with pattern from \cite{Liu2013}, and we assume that we only have access to an occluded image. 
Moreover, the image is contaminated by salt and pepper noise of density ${1}/{10}$. 
We seek to approximate the original from this noisy image. 

This is essentially a matrix completion problem, and  most of the scalable techniques rely on the Gaussian noise model. 
Note however the corresponding least-squares formulation is a good model against outliers:
\begin{equation*}
\label{eqn:example-matrix-completion}
\min_{X \in \mathcal{X}} ~\frac{1}{2}\norm{{A}(X)-b}^2 \quad \subjto \quad 0 \leq X \leq 1,  
\end{equation*}
where $\mathcal{X} = \{ X \in \R^{n \times n} : \norm{X}_{S_1} \leq \rho \}$ is a scaled nuclear norm ball, and $A:\R^{n \times n} \to \R^d$ is the sampling operator. 

Our framework also covers the following least absolute deviations formulation which is known to be more robust:
\begin{equation*}
\label{eqn:example-robust-matrix-completion}
\min_{X \in \mathcal{X}} ~\norm{{A}(X)-b}_1 \quad \subjto \quad 0 \leq X \leq 1. 
\end{equation*}
We solve both formulations with our framework, starting from all zero matrix, running $1000$ iterations, and assuming that we know the true nuclear norm of the original image. 
We choose $\beta_0 = 1$ in both cases. 

\begin{figure}[h]
\includegraphics[width=0.925\linewidth]{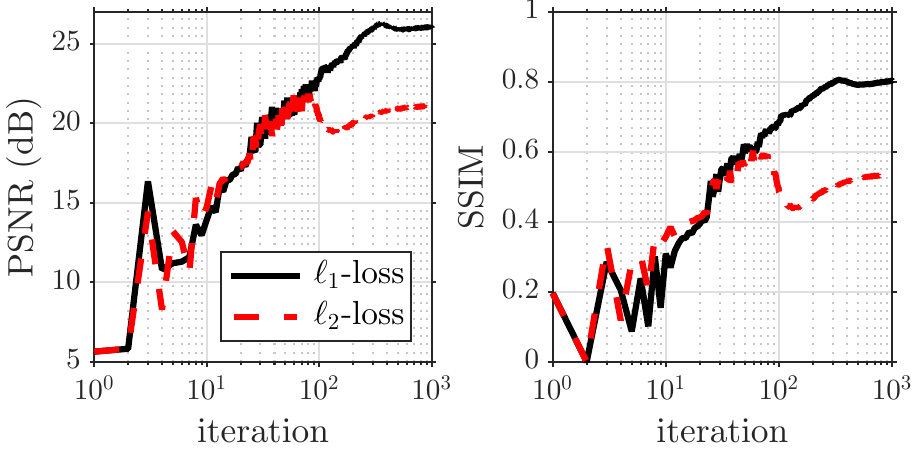}
\vspace{-0.25em}
        \caption{PSNR and SSIM vs iteration counter for formulations with $\ell_1$ and $\ell_2$ loss.}
        \label{fig:psnr-ssim}
        \vspace{-0.35em}
\end{figure}

This experiment demonstrates the implications of the flexibility of our framework in a simple machine learning setup. 
We compile the results in \Cref{fig:inpainting}, where the non-smooth formulation recovers a better approximation with $5$dB higher peak signal to noise ratio (PSNR) and $0.27$ higher structural similarity index (SSIM). 
Evaluation of PSNR and SSIM vs iteration counter are shown in \Cref{fig:psnr-ssim}. 

\section{Conclusions}\label{sec:conclusions}

We presented a CGM framework for the composite convex minimization template, that provably achieves the optimal rate. 
This rate also holds under approximate oracle calls with additive or multiplicative errors. 

Apart from its generalizations for various templates, there has been many attempts to improve the convergence rate, the arithmetic and the storage cost, or the proof techniques of CGM under some specific settings, \textit{cf.} \cite{Dunn1979, Guelat1986, Beck2004, Garber2015, JulienLacoste2015, Odor2016, Freund2016, Yurtsever2016} and references therein. 

Many of these techniques can be adapted to our framework, since we preserve key features of CGM, such as the reduced costs and atomic representations. 
The only seeming drawback is the loss of affine invariance, left for future, which is fundamentally challenging due to smoothing technique.

\section*{Acknowledgements} 
$^1$This work was supported by the Swiss National Science Foundation (SNSF) under grant number $200021\_178865 / 1$. 
$^1$This project has received funding from the European Research Council (ERC) under the European Union's Horizon $2020$ research and innovation programme (grant agreement no $725594$ - time-data). $^2$This work was supported by a public grant as part of the Investissement d'avenir project, reference ANR-$11$-LABX-$0056$-LMH, LabEx LMH, in a joint call with PGMO. $^{3~4}$This project has received funding from the Max Planck ETH Center for Learning Systems and by ETH core funding (to Gunnar R\"atsch).

{\small
\bibliography{../common/bibliography.bib}
\bibliographystyle{icml2018}
}

\newpage 
\onecolumn

\section*{Appendix}

\section*{A1~~~Preliminaries}

The following properties of smoothing are key to derive the convergence rate of our algorithm. 

Let $g:\R^d \to \R \cap \{+\infty\}$ be a a proper, closed and convex function, and denote its smooth approximation by
\begin{equation*}
g_{\beta} (z) = \max_{y \in \R^d} \ip{z}{y} - g^\ast(y) - \frac{\beta}{2} \norm{y}^2
\end{equation*}
where $g^\ast$ represents the Fenchel conjugate of $g$ and $\beta > 0$ is the smoothing parameter. 
Then, $g_\beta$ is convex and $\tfrac{1}{\beta}$-smooth. 
Let us denote the unique maximizer of this concave problem by 
\begin{align*}
y^\ast_\beta(z) & = \arg\max_{y \in \R^d} \ip{z}{y} - g^\ast(y) - \frac{\beta}{2} \norm{y}^2 \\
& = \arg\min_{y \in \R^d} \frac{1}{\beta} g^\ast(y) - \frac{1}{\beta} \ip{z}{y} + \frac{1}{2} \norm{y}^2 + \frac{1}{2} \norm{\frac{1}{\beta} z}^2 \\
& = \arg\min_{y \in \R^d} \frac{1}{\beta} g^\ast(y) + \frac{1}{2} \norm{y - \frac{1}{\beta} z}^2 \\
& = \prox{\beta^{-1} g^\ast} (\beta^{-1} z) 
 = \frac{1}{\beta} \big( z - \prox{\beta g} (z) \big)
\end{align*}
where the last equality is known as the Moreau decomposition. 
Then, the followings hold for $\forall z_1, z_2 \in \R^d$ and $\forall \beta, \gamma > 0$
\begin{align}
g_\beta(z_1) & \geq g_\beta(z_2) + \ip{\nabla g_{\beta}(z_2)}{z_1 - z_2} + \frac{\beta}{2} \norm{y^\ast_\beta(z_2) - y^\ast_\beta(z_1)}^2 \label{eqn:smoothing-prop-1} \\
g(z_1) & \geq g_\beta(z_2) + \ip{\nabla g_{\beta}(z_2)}{z_1 - z_2} + \frac{\beta}{2} \norm{y^\ast_\beta(z_2)}^2 \label{eqn:smoothing-prop-2} \\
g_{\beta} (z_1) & \leq g_{\gamma}(z_1) + \frac{\gamma - \beta}{2} \norm{y^\ast_\beta(z_1)}^2 \label{eqn:smoothing-prop-3}
\end{align}
Proofs can be found in Lemma~10 from \cite{TranDinh2017}.

Suppose that $g$ is $L_g$-Lipschitz continuous. Then, for any $\beta > 0$ and any $z \in \R^d$, the following bound holds:
\begin{equation}\label{eqn:smoothing-sandwich}
g_{\beta} (z) \leq g (z) \leq g_{\beta} (z) + \frac{\beta}{2} L_g^2
\end{equation}
Proof follows from equation~(2.7) in \cite{Nesterov2005} with a remark on the duality between Lipshitzness and bounded support (\textit{cf.} Lemma~5 in \cite{Dunner2016}).  

\newpage

\section*{A2~~~Convergence analysis}

This section presents the proof of our convergence results. 
We skip proofs of \Cref{thm:composite,thm:lipschitz,thm:indicator} since we can get these results as a special case by setting $\delta = 0$ in \Cref{thm:composite-inexact,thm:lipschitz-inexact,thm:indicator-inexact}.

\subsection*{Proof of \Cref{thm:composite-inexact}}

First, we use the smoothness of $F_{\beta_k}$ to upper bound the progress. 
Note that $F_{\beta_k}$ is $(L_f + \norm{A}^2/\beta_k)$-smooth. 
\begin{align}\label{eq_proof:upper_bound_smooth}
F_{\beta_k}(x_{k+1}) &\leq F_{\beta_k}(x_k)  +\eta_k 
\ip{\nabla F_{\beta_k}(x_k) }{\tilde{s}_k - x_k}
+ \frac{\eta_k^2}{2}\norm{\tilde{s}_k - x_k}^2(L_f + \frac{\norm{A}^2}{\beta_k}) \nonumber \\
&\leq F_{\beta_k}(x_k)  +\eta_k 
\ip{\nabla F_{\beta_k}(x_k) }{\tilde{s}_k - x_k}
+ \frac{\eta_k^2}{2}D_{\mathcal{X}}^2(L_f + \frac{\norm{A}^2}{\beta_k}),
\end{align}
where $\tilde{s}_k$ denotes the atom selected by the inexact $\mathrm{lmo}$, and the second inequality follows since $\tilde{s}_k\in\mathcal{X}$.  

By definition of inexact oracle \eqref{eqn:inexact_LMO}, we have
\begin{align*}
\ip{ \nabla F_{\beta_k}(x_k)}{\tilde{s}_k - x_k}
&\leq \ip{\nabla F_{\beta_k} (x_k)}{s_k - x_k}
+ \delta \frac{\eta_k}{2}D_{\mathcal{X}}^2(L_f + \frac{\norm{A}^2}{\beta_k}) \\
&\leq \ip{\nabla F_{\beta_k} (x_k)}{x^\star - x_k}
+ \delta \frac{\eta_k}{2}D_{\mathcal{X}}^2(L_f + \frac{\norm{A}^2}{\beta_k}) \\
&= 
\ip{\nabla f(x_k)}{x^\star - x_k} + \ip{A^\top \nabla g_{\beta_k}(Ax_k)}{x^\star - x_k }
+ \delta \frac{\eta_k}{2}D_{\mathcal{X}}^2(L_f + \frac{\norm{A}^2}{\beta_k}) ,
\end{align*}
where the second line follows since $s_k$ is a solution of $\min_{x \in \mathcal{X}} \ip{\nabla F_{\beta_k} (x_k)}{x}$.

Now, convexity of $f$ ensures $\langle\nabla f(x_k),x^\star - x_k \rangle\leq f(x^\star) - f(x_k)$. 
Using property \eqref{eqn:smoothing-prop-2}, we have 
\begin{align*}
\ip{A^\top \nabla g_{\beta_k} (Ax_k)}{x^\star - x_k} 
&= 
\ip{\nabla g_{\beta_k} (Ax_k)}{Ax^\star - Ax_k} \\
&\leq 
g(Ax^\star) - g_{\beta_k}(Ax_k) - \frac{\beta_k}{2} \norm{  y^\ast_{\beta_k}(Ax_k) }^2. 
\end{align*}
Putting these altogether, we get the following bound
\begin{align}
F_{\beta_k}(x_{k+1}) 
&\leq 
F_{\beta_k}(x_k) +\eta_k \left( f(x^\star) - f(x_k) + g(Ax^\star) - g_{\beta_k}(Ax_k) - \frac{\beta_k}{2} \norm{ \nabla y^\ast_{\beta_k}(Ax_k) }^2 \right) \nonumber \\
& \qquad + \frac{\eta_k^2}{2}D_{\mathcal{X}}^2(L_f + \frac{\norm{A}^2}{\beta_k})(1+\delta) \label{eqn:proof-recursion-1} \\
& = 
(1 - \eta_k) F_{\beta_k}(x_k) +\eta_k  F(x^\star) - \frac{\eta_k \beta_k}{2} \norm{ \nabla y^\ast_{\beta_k}(Ax_k) }^2 
+ \frac{\eta_k^2}{2}D_{\mathcal{X}}^2(L_f + \frac{\norm{A}^2}{\beta_k})(1+\delta). \nonumber
\end{align}
Now, using \eqref{eqn:smoothing-prop-3}, we get
\begin{align*}
F_{\beta_k}(x_k) 
& = 
f(x_k) + g_{\beta_k}(Ax_k) \\
& \leq 
f(x_k) + g_{\beta_{k-1}}(Ax_k) + \frac{\beta_{k-1}-\beta_k}{2}\norm{y^\ast_{\beta_k}(Ax_k)}^2 \\
& = 
F_{\beta_{k-1}}(x_k) + \frac{\beta_{k-1}-\beta_k}{2}\norm{y^\ast_{\beta_k}(Ax_k)}^2.
\end{align*}
We combine this with \eqref{eqn:proof-recursion-1} and subtract $F(x^\star)$ from both sides to get
\begin{align*}
F_{\beta_k}(x_{k+1}) - F(x^\star)
&\leq 
(1 - \eta_k) \big( F_{\beta_{k-1}}(x_k) -  F(x^\star) \big) + \frac{\eta_k^2}{2}D_{\mathcal{X}}^2(L_f + \frac{\norm{A}^2}{\beta_k})(1+\delta) \\ 
& \qquad +  \big( (1-\eta_k) (\beta_{k-1} - \beta_k) - \eta_k \beta_k \big) \frac{1}{2} \norm{ y^\ast_{\beta_k}(Ax_k) }^2 .
\end{align*}

Let us choose $\eta_k$ and $\beta_k$ in a way to vanish the last term. 
By choosing $\eta_k = \frac{2}{k+1}$ and $\beta_k = \frac{\beta_0}{\sqrt{k+1}}$ for $k \geq 1$ with some $\beta_0 > 0$, we get $(1-\eta_k)(\beta_{k-1}-\beta_k) - \eta_k\beta_k < 0$. Hence, we end up with 
\begin{align*}
F_{\beta_k}(x_{k+1}) - F(x^\star)
\leq 
(1 - \eta_k) \big( F_{\beta_{k-1}}(x_k) -  F(x^\star) \big) + \frac{\eta_k^2}{2}D_{\mathcal{X}}^2(L_f + \frac{\norm{A}^2}{\beta_k})(1+\delta).
\end{align*}
By recursively applying this inequality, we get 
\begin{align*}
F_{\beta_k}(x_{k+1}) - F(x^\star) 
& \leq 
\prod_{j = 1}^{k} (1-\eta_j)\left( F_{\beta_{j-1}}(x_k) - F(x^\star) \right) + \frac{1}{2}D_{\mathcal{X}}^2(1+\delta) \sum_{\ell=1}^k \eta_{\ell}^2 \prod_{j=\ell}^k (1-\eta_j) (L_f + \frac{\|A\|^2}{\beta_\ell}) \\
& \leq 
\prod_{j = 1}^{k} (1-\eta_j)\left( F_{\beta_{j-1}}(x_k) - F(x^\star) \right) + \frac{1}{2}D_{\mathcal{X}}^2 (L_f + \frac{\|A\|^2}{\beta_k}) (1+\delta) \sum_{\ell=1}^k \eta_{\ell}^2 \prod_{j=\ell}^k (1-\eta_j) \\
& = 
\frac{1}{2}D_{\mathcal{X}}^2 (L_f + \frac{\|A\|^2}{\beta_k}) (1+\delta) \sum_{\ell=1}^k \eta_{\ell}^2 \prod_{j=\ell}^k (1-\eta_j),
\end{align*}
where the second line follows since $\beta_k \leq \beta_j$ for any positive integer $j \leq k$, and the third line since $\eta_1 = 1$.

Now, we use the following relation
\begin{equation*}
\sum_{\ell=1}^k \eta_{\ell}^2 \prod_{j=\ell}^k (1-\eta_j) = \sum_{\ell=1}^k \frac{4}{(\ell+1)^2} \prod_{j=\ell}^k \frac{j-1}{j+1} = \sum_{\ell=1}^k \frac{4}{(\ell+1)^2} \frac{(\ell-1)\ell}{k(k+1)} \leq \frac{4}{k+1},
\end{equation*}
which yields the first result of \Cref{thm:composite-inexact} as
\begin{align*}
F_{\beta_k}(x_{k+1}) - F(x^\star) 
& \leq 
\frac{2}{k+1}D_{\mathcal{X}}^2 (L_f + \frac{\|A\|^2}{\beta_k}) (1+\delta) = 
2D_{\mathcal{X}}^2 (\frac{L_f}{k+1} + \frac{\|A\|^2 }{\beta_0 \sqrt{k+1}}) (1+\delta).
\end{align*}

\subsection*{Proof of \Cref{thm:lipschitz-inexact}}

Now, we further assume that $g:\R^d \to \R \cup \{+\infty\}$ is $L_g$-Lipschitz continuous. 
From \eqref{eqn:smoothing-sandwich}, we get
\begin{equation*}
g(Ax_{k+1}) \leq g_{\beta_k}(Ax_{k+1}) + \frac{\beta_kL_g^2 }{2} = g_{\beta_k}(Ax_{k+1}) + \frac{\beta_0 L_g^2}{2\sqrt{k+1}} . \end{equation*}

We complete the proof by adding $f(x_{k+1}) - F(x^\star)$ to both sides:
\begin{align*}
F(x_{k+1}) - F(x^\star) 
\leq 
F_{\beta_k}(x_{k+1}) - F(x^\star) + \frac{\beta_0 L_g^2}{2\sqrt{k+1}}.
\end{align*}

\subsection*{Proof of \Cref{thm:indicator-inexact}}

From the Lagrange saddle point theory, we know that the following bound holds $\forall x \in \mathcal{X}$ and $\forall r \in \mathcal{K}$:
\begin{equation*}
f^\star \leq \mathcal{L}(x,r,y^\star) = f(x) + \ip{y_\star}{Ax - r} \leq f(x) + \norm{y_\star}\norm{Ax - r},
\end{equation*}
Since $x_{k+1} \in \mathcal{X}$, we get
\begin{equation}\label{eqn:obj-lower-bound}
f(x_{k+1}) - f^\star \geq - \min_{r\in\mathcal{K}}\norm{y^\star}\norm{Ax_{k+1} - r} =  - \norm{y^\star}\mathrm{dist}(Ax_{k+1}, \mathcal{K}). 
\end{equation}
This proves the first bound in \Cref{thm:indicator-inexact}.

The second bound directly follows by \Cref{thm:composite-inexact} as
\begin{align*}
f(x_{k+1}) - f^\star \leq \underbrace{f(x_{k+1}) - f^\star + \frac{1}{2\beta_k} \mathrm{dist}^2(Ax_{k+1}, \mathcal{K})}_{F_{\beta_k}(x_{k+1}) - F^\star} \leq 
2D_{\mathcal{X}}^2 (\frac{L_f}{k+1} + \frac{\|A\|^2 }{\beta_0 \sqrt{k+1}}) (1+\delta).
\end{align*}

Now, we combine this with \eqref{eqn:obj-lower-bound}, and we get
\begin{align*}
- \norm{y^\star}\mathrm{dist}(Ax_{k+1}, \mathcal{K}) + \frac{1}{2\beta_k} \mathrm{dist}^2(Ax_{k+1}, \mathcal{K}) 
& \leq 2D_{\mathcal{X}}^2 (\frac{L_f}{k+1} + \frac{\|A\|^2 }{\beta_0 \sqrt{k+1}}) (1+\delta) \\
& \leq 2D_{\mathcal{X}}^2 \frac{\beta_k}{\beta_0} (L_f + \frac{\|A\|^2 }{\beta_0}) (1+\delta) .
\end{align*}
This is a second order inequality in terms of $\mathrm{dist}(Ax_k, \mathcal{K})$. Solving this inequality, we get
\begin{align*}
\mathrm{dist}(Ax_{k+1}, \mathcal{K})  & \leq \beta_k \left( \norm{y^\star} + \sqrt{\norm{y^\star}^2 +  4 D_{\mathcal{X}}^2 \frac{1}{\beta_0} \left( L_f + \frac{\norm{A}^2}{\beta_0} \right)(1+\delta)}  \right) \\
& \leq \frac{2 \beta_0 }{\sqrt{k+1}} \left( \norm{y^\star} + D_{\mathcal{X}} \sqrt{\frac{1}{\beta_0}\left( L_f + \frac{\norm{A}^2}{\beta_0} \right)(1+\delta) } \right).
\end{align*}

\subsection*{Proof of \Cref{thm:composite-inexact-mult}}
Let us define the multiplicative error $\delta$ of the LMO:
\begin{align}\label{eqn:inexact_LMO2}
\ip{v_k}{ \tilde{s}_k - x_k} \leq \delta\ip{v_k}{s_k-x_k}
\end{align}
For the proof we assume that $x_1$ is feasible.
First, we use the smoothness of $F_{\beta_k}$ to upper bound the progress. 
Note that $F_{\beta_k}$ is $(L_f + \norm{A}^2/\beta_k)$-smooth. 
\begin{align}\label{eq_proof:upper_bound_smooth2}
F_{\beta_k}(x_{k+1}) &\leq F_{\beta_k}(x_k)  +\eta_k 
\ip{\nabla F_{\beta_k}(x_k) }{\tilde{s}_k - x_k}
+ \frac{\eta_k^2}{2}\norm{\tilde{s}_k - x_k}^2(L_f + \frac{\norm{A}^2}{\beta_k}) \nonumber \\
&\leq F_{\beta_k}(x_k)  +\eta_k 
\ip{\nabla F_{\beta_k}(x_k) }{\tilde{s}_k - x_k}
+ \frac{\eta_k^2}{2}D_{\mathcal{X}}^2(L_f + \frac{\norm{A}^2}{\beta_k}),
\end{align}
where $\tilde{s}_k$ denotes the atom selected by the inexact linear minimization oracle, and the second inequality follows since $\tilde{s}_k\in\mathcal{X}$.  

By definition of inexact oracle \eqref{eqn:inexact_LMO2}, we have
\begin{align*}
\ip{ \nabla F_{\beta_k}(x_k)}{\tilde{s}_k - x_k}
&\leq \delta \ip{\nabla F_{\beta_k} (x_k)}{s_k - x_k} \\
&\leq \delta\ip{\nabla F_{\beta_k} (x_k)}{x^\star - x_k}\\
&= 
\delta\ip{\nabla f(x_k)}{x^\star - x_k} + \delta\ip{A^\top \nabla g_{\beta_k}(Ax_k)}{x^\star - x_k },
\end{align*}
where the second line follows since $s_k$ is a solution of $\min_{x \in \mathcal{X}} \ip{\nabla F_{\beta_k} (x_k)}{x}$.

Now, convexity of $f$ ensures $\langle\nabla f(x_k),x^\star - x_k \rangle\leq f(x^\star) - f(x_k)$. 
Using property \eqref{eqn:smoothing-prop-2}, we have 
\begin{align*}
\ip{A^\top \nabla g_{\beta_k} (Ax_k)}{x^\star - x_k} 
&= 
\ip{\nabla g_{\beta_k} (Ax_k)}{Ax^\star - Ax_k} \\
&\leq 
g(Ax^\star) - g_{\beta_k}(Ax_k) - \frac{\beta_k}{2} \norm{  y^\ast_{\beta_k}(Ax_k) }^2. 
\end{align*}
Putting these altogether, we get the following bound
\begin{align}
F_{\beta_k}(x_{k+1}) 
&\leq 
F_{\beta_k}(x_k) +\eta_k \delta\left( f(x^\star) - f(x_k) + g(Ax^\star) - g_{\beta_k}(Ax_k) - \frac{\beta_k}{2} \norm{ \nabla y^\ast_{\beta_k}(Ax_k) }^2 \right) \nonumber \\
& \qquad\qquad + \frac{\eta_k^2}{2}D_{\mathcal{X}}^2(L_f + \frac{\norm{A}^2}{\beta_k})  \label{eqn:proof-recursion-1_2} \\
& = 
(1 - \delta\eta_k) F_{\beta_k}(x_k) +\delta\eta_k  F(x^\star) - \frac{\delta\eta_k \beta_k}{2} \norm{ \nabla y^\ast_{\beta_k}(Ax_k) }^2 
+ \frac{\eta_k^2}{2}D_{\mathcal{X}}^2(L_f + \frac{\norm{A}^2}{\beta_k}). \nonumber
\end{align}
Now, using \eqref{eqn:smoothing-prop-3}, we get
\begin{align*}
F_{\beta_k}(x_k) 
& = 
f(x_k) + g_{\beta_k}(Ax_k) \\
& \leq 
f(x_k) + g_{\beta_{k-1}}(Ax_k) + \frac{\beta_{k-1}-\beta_k}{2}\norm{y^\ast_{\beta_k}(Ax_k)}^2 \\
& = 
F_{\beta_{k-1}}(x_k) + \frac{\beta_{k-1}-\beta_k}{2}\norm{y^\ast_{\beta_k}(Ax_k)}^2.
\end{align*}
We combine this with \eqref{eqn:proof-recursion-1_2} and subtract $F(x^\star)$ from both sides to get
\begin{align*}
F_{\beta_k}(x_{k+1}) - F(x^\star)
&\leq 
(1 - \delta\eta_k) \big( F_{\beta_{k-1}}(x_k) -  F(x^\star) \big) + \frac{\eta_k^2}{2}D_{\mathcal{X}}^2(L_f + \frac{\norm{A}^2}{\beta_k}) \\ 
& \qquad +  \big( (1-\delta\eta_k) (\beta_{k-1} - \beta_k) - \delta\eta_k \beta_k \big) \frac{1}{2} \norm{ \nabla y^\ast_{\beta_k}(Ax_k) }^2 .
\end{align*}

By choosing $\eta_k = \frac{2}{\delta (k-1)+2}$ and $\beta_k = \frac{\beta_0}{\sqrt{\delta k+1}}$ for some $\beta_0 > 0$, we get $(1-\delta\eta_k)(\beta_{k-1}-\beta_k) - \delta\eta_k\beta_k < 0$ for any $k \geq 1$, hence we end up with 
\begin{align*}
F_{\beta_k}(x_{k+1}) - F(x^\star)
\leq 
(1 - \delta\eta_k) \big( F_{\beta_{k-1}}(x_k) -  F(x^\star) \big) + \frac{\eta_k^2}{2}D_{\mathcal{X}}^2(L_f + \frac{\norm{A}^2}{\beta_k}).
\end{align*}
Let us call for simplicity $C := D_{\mathcal{X}}^2(L_f + \frac{\norm{A}^2}{\beta_k})$, $E_{k+1} := F_{\beta_k}(x_{k+1}) - F(x^\star)$
Therefore, we have
\begin{align}\label{eqn:induction}
E_{k+1} \leq (1 - \delta\eta_k)E_k + \frac{\eta_k^2}{2}C
\end{align}

We now show by induction that:
\begin{align*}
E_{k} \leq 2\frac{\frac{1}{\delta}C+E_1}{\delta(k-1)+2}
\end{align*} 
The base case $k=1$ is trivial as $C>0$.
Call for simplicity $K := \delta(k-1)+2$. Note that $K\geq 2$. Under this notation we can write $\eta_k = \frac{2}{\delta (k-1)+2} = \frac{2}{K}$
For the induction step, we add a positive term ($E_1$ is positive as $x_1$ is assumed feasible) to \eqref{eqn:induction} and use the induction hypothesis:
\begin{align*}
E_{k+1} &\leq (1 - \delta\eta_k)E_k + \frac{\eta_k^2}{2}C + 2\delta\frac{E_1}{K^2}\\
&\leq (1 - \delta\frac{2}{K})E_k + \frac{2}{K^2}C + 2\delta\frac{E_1}{K^2}\\
&\leq (1 - \delta\frac{2}{K})2\frac{\frac{1}{\delta}C+E_1}{K} + \frac{2}{K^2}C + 2\delta\frac{E_1}{K^2}\\
&= (1 - \delta\frac{2}{K})2\frac{\frac{1}{\delta}C+E_1}{K} + 2\delta\left(\frac{\frac{1}{\delta}C}{K^2} + \frac{E_1}{K^2}\right)\\
&= 2\frac{\frac{1}{\delta}C+E_1}{K} \left(1 - \delta\frac{2}{K} + \frac{\delta}{K}\right)\\
&= 2\frac{\frac{1}{\delta}C+E_1}{K} \left(1 - \frac{\delta}{K}\right)\\
& \leq 2\frac{\frac{1}{\delta}C+E_1}{K+\delta}
\end{align*}
noting that $K+\delta = \delta k+2$ concludes the proof.

Proof of \Cref{thm:lipschitz-inexact-mult,thm:indicator-inexact-mult} follows similarly to the proofs of \Cref{thm:lipschitz-inexact,thm:indicator-inexact}. 

\newpage

\section*{A3~~~Numerical illustration for the pathological example in \Cref{sec:non-smooth}}

This section considers the pathological example by \citet{Nesterov2017} that we present in \Cref{sec:non-smooth}. 
We numerically demonstrate that our framework successfully finds the solution in contrast to the classical CGM. 

\begin{figure}[ht]
\centering
\includegraphics[width=0.65\linewidth]{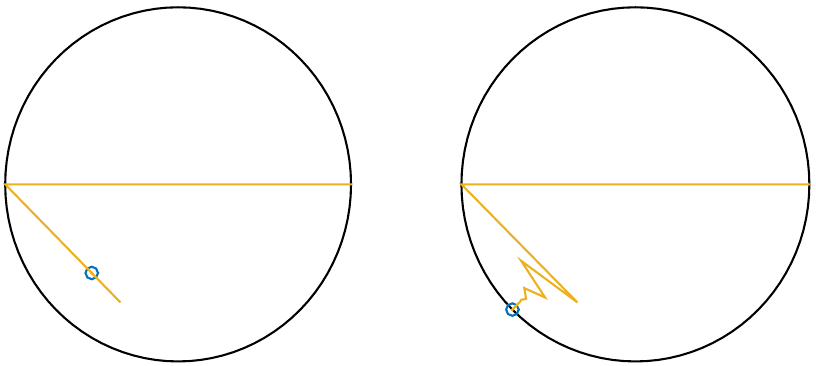}
        \caption{Classical CGM (\textit{left}) vs our framework (\textit{right}) for the pathological example, starting from $\bfrac{1}{0}$.}
        \label{fig:concept_pr}
\end{figure}

\Cref{fig:concept_pr} illustrates the paths followed by classical CGM vs our framework, starting both methods from $[1, 0]^\top$. Recall that the unique solution is $x^\star = [-\frac{1}{\sqrt{2}}, -\frac{1}{\sqrt{2}}]^\top$, but classical CGM converges to $[-\frac{1}{2}, -\frac{1}{2}]^\top$, which belongs to the convex hull of $\{x_0, \bfrac{\text{-}1}{0}, \bfrac{0}{\text{-}1} \}$. 
Our framework, on the other hand, converges to $x^\star$. 

\begin{figure}[ht]
\centering
\includegraphics[width=0.55\linewidth]{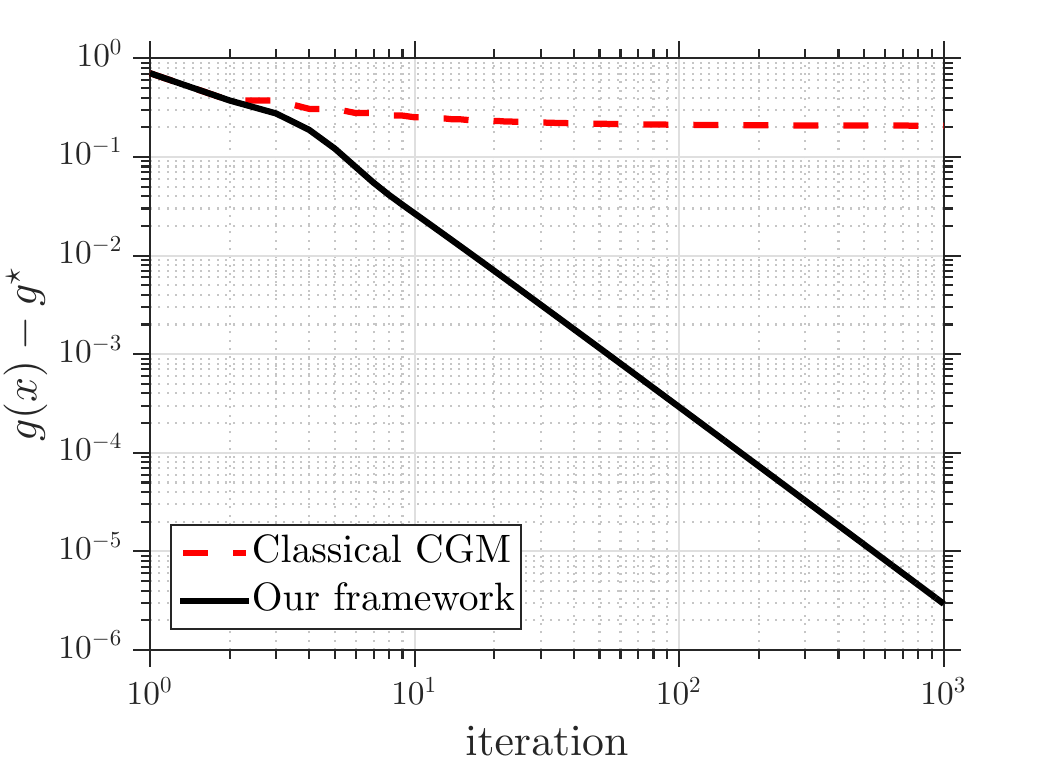}
        \caption{Evolution of the objective residual as a function of the iteration counter for the pathological example, starting from $\bfrac{1}{0}$.}
        \label{fig:concept_pr_convergence}
\end{figure}

Finally, \Cref{fig:concept_pr_convergence} plots the evolution of the objective residual as a function of the iteration counter. 
We observe an empirical $\mathcal{O}(1/k^2)$ convergence rate in this particular instance. 

\end{document}